\documentclass[11pt]{article}

 \newtheorem{theorem}{Theorem}
  \newtheorem{lemma}{Lemma}
 
 \newtheorem{corollary}{Corollary}
 
 \newtheorem{algorithm}{Algorithm}
\begin{document}
\title{Obtaining Exact Value by Approximate Computations
\thanks{The work is partially supported by China 973 Project NKBRPC-2004CB318003. Article was submitted to Science in China in June, 2006}}
\author{
Jingzhong  Zhang \ \ \ \  Yong Feng\\
Laboratory for Automated Reasoning and Programming\\
Chengdu Institute of Computer Applications\\
Chinese Academy of Sciences\\
610041 Chengdu, P. R. China\\
E-mail: zjz101@yahoo.com.cn, yongfeng@casit.ac.cn }
\date{}
\maketitle

\begin{abstract}
Numerical approximate computation can solve large and complex
problems fast. It has the advantage of high efficiency. However it
only gives approximate results, whereas  we need exact results in
many fields. There is a gap between approximate computation and
exact results. In this paper, we build a bridge by which exact
results can be obtained by numerical approximate computations.

\end{abstract}

\section{Introduction}
Historically, the fields of symbolic computation and numerical
computation have been developed by two distinct groups of people,
having relatively little interaction and overlaps with each other.
Symbolic computations are principally exact and stable. However
they have a high complexity. Therefore, they are slow and in
practice, are applicable only to small systems. Numerical
approximate computation has the advantage of being fast, flexible
in accuracy and being applicable to large scale problems. In
Recent two decades, numerical methods are applied in the field of
symbolic computations. In 1985, Kaltofen presented an algorithm
for performing the absolute irreducible factorization, and
suggested to perform his algorithm by floating-point numbers, then
the factor obtained is an approximate one. After then, numerical
methods have been studied to get approximate factors of a
polynomial\cite{corless_C}\cite{Huang}\cite{Mou}\cite{Sasaki_a}
\cite{Sasaki_b}\cite{Sasaki_c}. In the meantime, numerical methods
are applied to get approximate greatest common divisors of
approximate polynomials
\cite{beckermann}\cite{Corless_A}\cite{karmarkar}\cite{corless_D},
to compute functional decompositions\cite{corless_E}, to test
primality\cite{Galligo} and to find zeroes of a
polynomial\cite{Reid}. In 2000, Corless et al. applied numerical
method in implicitization of parametric curves, surfaces and
hypersurfaces\cite{corless_B}. The resulting implicit equation is
still an approximate one.

There is a gap between approximate computations and exact
results\cite{yang}. People usually use rational number
computations to override the gap\cite{Feng}. In fact, these are
not approximate computations but big number computations, which
are also exact computations. In 2005, Zhang et al proposed an
algorithm to get exact factors of a multivariate polynomial by
approximate computation\cite{zhang} but they did not discuss how
to override the gap. Command {\it convert} in maple can obtain an
approximate rational number from a float if we set variable {\it
Digits} to a positive integer. However, in order to obtain exact
rational number from its approximation, we need to know two
things. One is at what accuracy the float should be obtained by
numerical method; another is when we should stop and return the
rational number we want by the continued fraction method. So we
can not obtain the exact rational number from its approximation by
command {\it convert} and variable {\it Digits}. In this paper, we
solve the two problems, which can be described as follows:

There is an unknown rational number $m/n$  we want to obtain, and
assume that there are approximate methods to obtain its
approximation at arbitrary accuracy. We also know an upper bound
$N$ of absolute value of its denominator in advance. The two
problems will be solved such as: At first, we discuss how to
determine $\varepsilon$ which is a function in $N$. And then use
the approximate methods to obtain a floating-point number $x$, an
approximation of $m/n$ at accuracy $\varepsilon>0$, i.e.
$|x-m/n|<\varepsilon$. Second, we give a criteria to stop our
program and return the exact rational number we want by continued
fraction method.

The remainder of the paper is organized as follows. Section 2
gives a review of continued fraction. Section 3 discusses how
small the error needs to ensure the exact number to be obtained
and how to get the exact number from its approximation. Section 4
gives some experimental results. The final section makes
conclusions.

\section{Continued fraction}
A continued fraction representation of a real number is one of the
forms:
\begin{equation}\label{equ:1}
 a_0+\frac{1}{a_1+\frac{1}{a_2+\frac{1}{a_3+\cdots}}},
\end{equation}
where $a_0$ is an integer and $a_1,a_2,a_3,\cdots$ are positive
integers. One can abbreviate the above continued fraction as
$[a_0;a_1,a_2,\cdots]$. For finite continued fractions, note that
$$[a_0;a_1,a_2,\cdots,a_n,1]=[a_0;a_1,a_2,\cdots,a_n+1].$$
So, for every finite continued fraction, there is another finite
continued fraction that represents the same number. Every finite
continued fraction is rational number and every rational number
can be represented in precisely two different ways as a finite
continued fraction. The other representation is one element
shorter, and the final term must be greater than 1 unless there is
only one element. However, every infinite continued fraction is
irrational, and every irrational number can be represented in
precisely one way as an infinite continued fraction. An infinite
continued fraction representation for an irrational number is
mainly useful because its initial segments provide excellent
rational approximations to the number. These rational numbers are
called the convergents of the continued fraction. Even-numbered
convergents are smaller than the original number, while
odd-numbered ones are bigger. If successive convergents are found,
with numerators $h_1,h_2,\cdots$, and denominators
$k_1,k_2,\cdots$, then the relevant recursive relation is:
$$h_n=a_nh_{n-1}+h_{n-2},\,\, k_n=a_nk_{n-1}+k_{n-2}.$$
The successive convergents are given by the formula
$$\frac{h_n}{k_n}=\frac{a_nh_{n-1}+h_{n-2}}{a_nk_{n-1}+k_{n-2}},$$
where $h_{-1}=1$, $h_{-2}=0$, $k_{-1}=0$ and $k_{-2}=1$. Here are
some useful theorems\cite{Fraction}:
\begin{theorem}
For any positive $x\in R$, it holds that
\begin{equation}
[a_0,a_1,\cdots,a_{n-1},x]=\frac{xh_{n-1}+h_{n-2}}{xk_{n-1}+k_{n-1}}
\end{equation}
\end{theorem}
\begin{theorem} \label{theo:2}
The convergents of $[a_0,a_1,a_2,\cdots]$ are given by
$$[a_0,a_1,\cdots,a_n]=\frac{h_n}{k_n}$$
and $$k_nh_{n-1}-k_{n-1}h_n=(-1)^n.$$
\end{theorem}
\begin{theorem}
$$\left |\frac{h_n}{k_n}-\frac{h_{n-1}}{k_{n-1}}\right |=\frac{1}{k_nk_{n-1}}$$
and
$$\frac{1}{k_n(k_{n+1}+k_n)}<\left |x-\frac{h_n}{k_n}\right|<\frac{1}{k_nk_{n+1}}$$
\end{theorem}

In order to recover exact rational number, we introduce a
controlling error into the conventional continued fraction method.
The continued fraction method is modified as follows.
\newcounter{num}
\begin{algorithm}\label{algor:1} Continued fraction method \\
Input: a nonnegative floating-point number $a$ and $\varepsilon>0$; \\
Output: a rational number $b$.
\begin{list}{Step \arabic{num}:}{\usecounter{num}\setlength{\rightmargin}{\leftmargin}}
\item $i:=1$ and $x_1:=a$; \item \label{step1:2}Getting integral
part of $x_i$ and assigning it to $a_i$, assigning its remains to
$b_i$. If $b_i<\varepsilon$, then goto Step \ref{step1:5}; \item
$i:=i+1$; \item $x_i:=\frac{1}{b_{i-1}}$ and goto Step
\ref{step1:2}; \item \label{step1:5}  Computing expression
(\ref{equ:1}) and assigning it to $b$. \item return $b$.
\end{list}
\end{algorithm}
We will discuss the controlling error $\varepsilon$ in algorithm
\ref{algor:1} in the next section.

\section{Recovering the exact number from its approximation}
In this section, we will solve such a problem: for a given
floating number $w$ which is an approximation of rational number
$\frac{m}{n}$, how do we obtain integer $m$ and $n$?  Without loss
of generality, we always assume that $m$,$n$ are positive number.
At first, we have a lemma as follows:
\begin{lemma} \label{lemma:1}
$m,n,p,q$ are integer, and $pn>0$. If
$|\frac{m}{n}-\frac{q}{p}|<\frac{1}{pn}$, then
$\frac{m}{n}=\frac{q}{p}$.
\end{lemma}
{\bf Proof}. $|\frac{m}{n}-\frac{q}{p}|=\frac{|pm-qn|}{pn}$.
Noticing $|pm-qn|$ is a nonnegative integer, and
$|\frac{m}{n}-\frac{q}{p}|<\frac{1}{pn}$, yields $|pm-qn|<1$.
Hence $|pm-qn|=0$. That is $\frac{m}{n}=\frac{q}{p}$. The proof is
finished.
\begin{corollary}\label{cor:1}
$m,n,p,q$ are integers and $p>0$,$n>0$. Let $N\ge\max\{p,n,2\}$.
If $|\frac{m}{n}-\frac{q}{p}|<\frac{1}{N(N-1)}$, then
$\frac{m}{n}=\frac{q}{p}$.
\end{corollary}
{\bf Proof}: When $p\ne n$, it holds that $pn\le N(N-1)$. Hence
$|\frac{m}{n}-\frac{q}{p}|<\frac{1}{N(N-1)}\le \frac{1}{pn}$.
According to lemma \ref{lemma:1}, it is obtained that
$\frac{m}{n}=\frac{q}{p}$. When $p=n$, we have
\begin{eqnarray*}
&&|\frac{m}{n}-\frac{q}{p}|=\frac{|m-q|}{n}<\frac{1}{N(N-1)}\Rightarrow
|m-q|<\frac{n}{N(N-1)}\le \frac{1}{N-1}\le 1
\end{eqnarray*}
So, $m=q$. The proof of the corollary is finished.

And now, we study how small the error $|w-\frac{m}{n}|$ needs so
as to get exact rational number $\frac{m}{n}$ from its
approximation $w$ . The following theorem answers this question.
\begin{theorem}\label{theo:add1}
Let $x=\frac{m}{n}$ be a reduced proper fraction, and
$N\ge\max\{n,2\}$. Assume that $|x-w|<1/(2N(N-1))$. If we get
positive rational number $p/q$ such that $|p/q-w|<1/(2N(N-1))$ ,
where $q\le N$ , then it holds that $x=q/p$.
\end{theorem}
{\bf Proof}: From the assumption of the theorem, we have
$|x-q/p|<1/(N(N-1))$.  According to corollary \ref{cor:1}, it
holds that $q/p=m/n=x$. The proof of the theorem is finished.

We have answered how small the error is so that we can recover the
exact rational number from its approximation. The remaining
problem is how to get the exact number. We attack it by  continued
fraction method.

Let $n_2/n_1$ be a rational number and $r_0$ its approximation.
Their continued fraction representations are
$n_2/n_1=[a_0,a_1,\cdots,a_L]$ and $r_0=[b_0,b_1,\cdots,b_M]$
respectively. We wish that $a_i=b_i$ for $i=0,1,2,\cdots,L-1$ and
for the last term of the continued fraction representations of
$n_2/n_1$, either $a_L=b_L$ or $a_L-1=b_L$, so that we can get
$n_2/n_1$ from $[b_0,b_1,\cdots,b_{L+1}]$. This is the following
theorem:
\begin{theorem}\label{theo:3}
Let $n_2/n_1$ be a rational number and $r_0$ its approximation.
Assume that $n_2$,$n_1$ are coprime positive numbers, where
$n_2<n_1$,and $n_1>1$. The representations of $n_2/n_1$ and $r_0$
are $[a_0,a_1,\cdots,a_L]$ and $[b_0,b_1,\cdots,b_M]$
respectively. If $|r_0-n_2/n_1|<1/(4n_1(n_1-1))$, then one of the
following statements must hold.
\begin{itemize}
\item $a_i=b_i$  ($i=0,1,\cdots,L$); \item $a_i=b_i$
($i=0,1,\cdots,L-1$),  $a_L-1=b_L$, and $b_{L+1}=1$.
\end{itemize}
\end{theorem}
According to assumption of $n_2<n_1$, we have that $a_0=0$, and
$b_0=0$. Hence $a_0=b_0$. In order to finish the  proof of theorem
\ref{theo:3}, we need two lemmas. Due to $
n_2/n_1=[a_0,a_1,\cdots,a_L]$ and $r_0=[b_0,b_1,\cdots,b_M]$, we
have the following expansions:
\begin{equation}\label{equ:add1}
\frac{n_1}{n_2}=a_1+\frac{n_3}{n_2},\;
\frac{n_2}{n_3}=a_2+\frac{n_4}{n_3},\; \cdots ,
\frac{n_{L-1}}{n_L}=a_{L-1}+\frac{1}{n_L},\; n_L=a_L
\end{equation}
and
\begin{equation} \label{equ:add2}
\frac{1}{r_0}=b_1+r_1,\; \frac{1}{r_1}=b_2+r_2,\; \cdots ,
\frac{1}{r_{L-1}}=b_L+r_L,\; \cdots,\; \frac{1}{r_{M-1}}=b_M
\end{equation}
Denoting $d_i=r_i-n_{i+2}/n_{i+1}$, we have a lemma as follows:
\begin{lemma} \label{lemma:2}
Let $n_2/n_1$ be a rational number and $r_0$ its approximation.
Assume that $n_2$,$n_1$ are coprime positive integers, where
$n_2<n_1$,and $n_1>1$. The representations of $n_2/n_1$ and $r_0$
are $[a_0,a_1,\cdots,a_L]$ and $[b_0,b_1,\cdots,b_M]$
respectively. And assume that $a_i=b_i$ for $i\le k<L$($k$ is a
positive integer). Then when $|d_k|<\frac{1}{n_{k+1}(n_{k+1}-1)}$,
it holds that $a_{k+1}=b_{k+1}$ for $k<L-1$; when
$|d_{L-1}|<\frac{1}{n_L(n_L+1)}$, it holds that $a_L=b_L$ or
$a_L-1=b_L$.
\end{lemma}
{\bf Proof}: At first, we show that under the assumption of the
lemma if we have
\begin{equation}\label{equ:add3}
\left|\frac{n_{k+1}^2d_k}{n_{k+2}(n_{k+2}+n_{k+1}d_k)}\right|<
\frac{1}{n_{k+2}}
\end{equation}
then, it holds that $a_{k+1}=b_{k+1}$ for $k<L-1$, and
$a_{k+1}=b_{k+1}$ or $a_{k+1}-1=b_{k+1}$ for $k=L-1$. We discuss
it in two cases:\\ Case 1($k<L-1$): From
$d_k=r_k-n_{k+2}/n_{k+1}$, it holds that
$r_k=n_{k+2}/n_{k+1}+d_k$. Hence we have that
\begin{eqnarray*}
&&\frac{1}{r_k}-\frac{n_{k+1}}{n_{k+2}}=-\frac{n_{k+1}^2d_k}{n_{k+2}(n_{k+2}+n_{k+1}d_k)}\Rightarrow
\frac{1}{r_k}=\frac{n_{k+1}}{n_{k+2}}-\frac{n_{k+1}^2d_k}{n_{k+2}(n_{k+2}+n_{k+1}d_k)}\\
&&\Rightarrow
\frac{1}{r_k}=a_{k+1}+\frac{n_{k+3}}{n_{k+2}}-\frac{n_{k+1}^2d_k}{n_{k+2}(n_{k+2}+n_{k+1}d_k)}=b_{k+1}+r_{k+1}
\end{eqnarray*}
Hence, it is obvious that $a_{k+1}=b_{k+1}$ if and only if
\begin{equation} \label{equ:2} 0\le
\frac{n_{k+3}}{n_{k+2}}-\frac{n_{k+1}^2d_k}{n_{k+2}(n_{k+2}+n_{k+1}d_k)}<1.
\end{equation}
Therefore, if inequality (\ref{equ:add3}) holds, then above
inequality is guaranteed.\\
Case 2:(when $k=L-1$) We have
$$\frac{1}{r_{L-1}}=a_L-\frac{n_{L}^2d_{L-1}}{n_{L+1}(n_{L+1}+n_{L}d_{L-1})}=b_L+r_L $$
From the above equation, if
$$|\frac{n_{L}^2d_{L-1}}{n_{L+1}(n_{L+1}+n_{L}d_{L-1})}|<\frac{1}{n_{L+2}}=1$$,
then $a_L=b_L$ for $d_{L-1}<0$, and $a_L-1=b_L$ for $d_{L-1}\ge
0$. Therefore, we have shown that if inequality (\ref{equ:add3})
holds, then $a_{k+1}=b_{k+1}$ for $k<L-1$, and $a_{k+1}=b_{k+1}$
or $a_{k+1}-1=b_{k+1}$ for $k=L-1$.

On the other hand, we have
\begin{eqnarray*}
&&\left|\frac{n_{k+1}^2d_k}{n_{k+2}(n_{k+2}+n_{k+1}d_k)}\right.|=\frac{n_{k+1}^2|d_k|}{|n_{k+2}(n_{k+2}+n_{k+1}d_k)|}\\
&&=
\frac{n_{k+1}|d_k|}{n_{k+2}(|\frac{n_{k+2}}{n_{k+1}}+d_k|)}\leq
\frac{n_{k+1}|d_k|}{n_{k+2}|(\frac{n_{k+2}}{n_{k+1}}-|d_k|)|}
\end{eqnarray*}
So, in order to ensure inequality (\ref{equ:add3}), we only need
it holds that
\begin{equation}\label{equ:3}
\frac{n_{k+1}|d_k|}{n_{k+2}|(\frac{n_{k+2}}{n_{k+1}}-|d_k|)|}<
\frac{1}{n_{k+2}}
\end{equation}
Solving inequality (\ref{equ:3}) yields
\begin{equation}\label{equ:4}
|d_k|<\frac{n_{k+2}}{n_{k+1}(n_{k+1}+1)}
\end{equation}
When $k<L-1$, we have that $n_{k+2}>1$. So, it holds that
$$\frac{1}{n_{k+1}(n_{k+1}-1)}\leq \frac{n_{k+2}}{n_{k+1}(n_{k+1}+1)} $$
Accordingly, it is obtained that
\begin{equation}
|d_k|<\frac{1}{n_{k+1}(n_{k+1}-1)}
\end{equation}
When $k=L-1$, we have that $n_{L+1}=1$, so it is obtained that
\begin{equation}
|d_{L-1}|<\frac{1}{n_L(n_L+1)}
\end{equation}

The proof of lemma \ref{lemma:2} is finished.

\begin{lemma} \label{lemma:relation_d}
Let $n_2/n_1$ be a rational number and $r_0$ its approximation,
where $n_2$,$n_1$ are coprime positive integers, and $n_2<n_1$,and
$n_1>1$. The continued fraction representations of $n_2/n_1$ and
$x$ are $[a_0,a_1,\cdots,a_L]$ and $[b_0,b_1,\cdots,b_M]$
respectively. Denote $d_i=r_i-n_{i+2}/n_{i+1}$ for $i=0,\cdots,L$.
Assume that $a_i=b_i$ for $i\le k<L-1$($k$ is a positive integer
). Then when $|d_k|<\frac{1}{n_{k+1}(n_{k+1}-1)}$, it holds that
\begin{equation}\label{equ:7}
|d_{k+1}|<\frac{n_{k+1}(n_{k+1}-1)}{n_{k+2}(n_{k+2}-1)}|d_k|
\end{equation}
\end{lemma}
{\bf Proof}: Under the assumption that $a_i=b_i$ for
$i=0,1,\cdots,k$, from equation (\ref{equ:add3}), we get
$d_{k+1}=-\frac{n_{k+1}^2d_k}{n_{k+2}(n_{k+2}+n_{k+1}d_k)}$. Hence
we deduce a relation as follows:
\begin{eqnarray*}
|d_{k+1}|&=&\frac{n_{k+1}^2|d_k|}{n_{k+2}^2+n_{k+1}n_{k+2}d_k}=\frac{n_{k+1}|d_k|}{n_{k+2}(\frac{n_{k+2}}{n_{k+1}}+d_k)}\\
&=&\frac{n_{k+1}(n_{k+1}-1)|d_k|}{n_{k+2}(\frac{n_{k+2}(n_{k+1}-1)}{n_{k+1}}+(n_{k+1}-1)d_k)}\\
&=&\frac{n_{k+1}(n_{k+1}-1)|d_k|}{n_{k+2}(n_{k+2}-1+\frac{n_{k+1}-n_{k+2}}{n_{k+1}}+(n_{k+1}-1)d_k)}
\end{eqnarray*}
When $|d_k|<\frac{1}{n_{k+1}(n_{k+1}-1)}$, it holds that
$\frac{n_{k+1}-n_{k+2}}{n_{k+1}}+(n_{k+1}-1)d_k>0$. Hence we have
a relation between $d_{k+1}$ and $d_k$:
\[
|d_{k+1}|<\frac{n_{k+1}(n_{k+1}-1)}{n_{k+2}(n_{k+2}-1)}|d_k| \]
The proof of the lemma is finished.

And now, let us prove the theorem. If
$|d_0|=|r_0-n_2/n_1|<1/(4n_1(n_1-1))$, From lemma
\ref{lemma:relation_d}, we can get $|d_i|<1/(4n_{i+1}(n_{i+1}-1))$
for $i=0,\cdots,L-1$. Note that $n_L>n_{L+1}=1$ and
$$\frac{1}{4n_{i+1}(n_{i+1}-1)}<\frac{1}{n_{i+1}(n_{i+1}+1)}<\frac{1}{n_{i+1}(n_{i+1}-1)}$$ when
$n_{i+1}>1$. So, it holds that
$$|d_i|<\frac{1}{4n_{i+1}(n_{i+1}-1)}<\frac{1}{n_{i+1}(n_{i+1}+1)}$$
for $i=0,\cdots,L-1$. According to lemma \ref{lemma:2}, the proof
of the theorem is finished.

For an unknown rational number $n_2/n_1$ and its approximation
$r_0$, theorem \ref{theo:3} shows that $n_2/n_1=[b_0,\cdots,b_L]$
or $n_2/n_1=[b_0,b_1,\cdots,b_L,1]$ when
$|r_0-n_2/n_1|<1/(4n_1(n_1-1))$. However, we do not know what the
number $L$ is. If we make $b_{L+1}$ large enough when
$n_2/n_1=[b_0,\cdots,b_L]$, or make $b_{L+2}$ large enough when
$n_2/n_1=[b_0,b_1,\cdots,b_L,1]$, then we recover $n_2/n_1$
easily. The following theorem solve this problem.
\begin{theorem}\label{theo:basis_alg}
Let $n_2/n_1$ be a rational number and $r_0$ its approximation.
Assume that $n_2$,$n_1$ are coprime positive integers, where
$n_2<n_1$,and $n_1>1$. $K$ is a positive integer. The continued
fraction representations of $n_2/n_1$ and $r_0$ are
$[a_0,a_1,\cdots,a_L]$ and $[b_0,b_1,\cdots,b_M]$ respectively. If
$|d_0|=|r_0-n_2/n_1|<1/((2K+2)n_1(n_1-1))$, then one of the
following two statements must hold
\begin{itemize}
\item $a_i=b_i$ for $i=0,\cdots,L$, and $b_{L+1}\ge K$; \item
$a_i=b_i$ for $i=0,\cdots,L-1$, and $b_L=a_L-1$, $b_{L+1}=1$,
$b_{L+2}\ge K$.
\end{itemize}
\end{theorem}
{\bf Proof}: From equation (\ref{equ:add1}) and equation
(\ref{equ:add2}), we have that
$\frac{n_{L-1}}{n_L}=a_{L-1}+\frac{1}{n_L},\; n_L=a_L$ and
$\frac{1}{r_{L-2}}=b_{L-1}+r_{L-1}$. When
$|d_0|=|r_0-n_2/n_1|<1/((2K+2)n_1(n_1-1))<1/(4n_1(n_1-1))$, from
lemma \ref{lemma:relation_d}, it holds that
$|d_i|<1/((2K+2)n_{i+1}(n_{i+1}-1))$ for $i=0,1,\cdots, L-1$.
Furthermore, from theorem \ref{theo:3} it holds that
$a_{L-1}=b_{L-1}$ and $d_{L-1}=r_{L-1}-\frac{1}{n_L}$. So we have
$$\frac{1}{r_{L-1}}=a_L-\frac{n_L^2d_{L-1}}{1+n_Ld_{L-1}}$$

We discuss it in two cases:\\
Case 1 ($d_{L-1}>0$): Due to
$-\frac{n_L^2d_{L-1}}{1+n_Ld_{L-1}}<0$, we have
$$\frac{1}{r_{L-1}}=a_L-1+(1-\frac{n_L^2d_{L-1}}{1+n_Ld_{L-1}})=b_L+r_L$$
Hence
\begin{eqnarray*}
&&r_L=(1-\frac{n_L^2d_{L-1}}{1+n_Ld_{L-1}})=\frac{1+n_Ld_{L-1}-n_L^2d_{L-1}}{1+n_Ld_{L-1}}\\
&&\Rightarrow
\frac{1}{r_L}=\frac{1+n_Ld_{L-1}}{1+n_Ld_{L-1}-n_L^2d_{L-1}}
\end{eqnarray*}
It is obvious that $1<1/r_L<2$ when $|d_{L-1}|<1/(4n_L(n_L-1))$.
Therefore, we have
$$\frac{1}{r_L}=1+\frac{1+n_Ld_{L-1}}{1+n_Ld_{L-1}-n_L^2d_{L-1}}-1=1+\frac{n_L^2d_{L-1}}{1+n_Ld_{L-1}-n_L^2d_{L-1}}$$
We get that $b_{L+1}=1$ and
\begin{equation}
\frac{1+n_Ld_{L-1}-n_L^2d_{L-1}}{n_L^2d_{L-1}}=b_{L+2}+r_{L+2}.
\end{equation}
So,
$$b_{L+2}=\left[\frac{1+n_Ld_{L-1}-n_L^2d_{L-1}}{n_L^2d_{L-1}}\right]
$$
where $[.]$ stands for getting the integral part of a number.

We hope that $b_{L+2}$ is greater than some integer $K$, which is
used as a sign that $b_0,\cdots,b_L,1$ have been obtained. Solve
the following inequality:
\begin{eqnarray*}
&&\frac{1+n_Ld_{L-1}-n_L^2d_{L-1}}{n_L^2d_{L-1}}>K\Leftrightarrow
1+n_Ld_{L-1}-n_L^2d_{L-1}>Kn_L^2d_{L-1} \\
&&\Leftrightarrow d_{L-1}<\frac{1}{(K+1)n_L^2-n_L}
\end{eqnarray*}
Since $(2K+2)n_L(n_L-1)>(K+1)n_L^2-n_L$, we take
$$d_{L-1}<\frac{1}{(2K+2)n_L(n_L-1)}$$

Case 2($d_{L-1}<0$): Due to
$-\frac{n_L^2d_{L-1}}{1+n_Ld_{L-1}}>0$, we have
$$\frac{1}{r_{L-1}}=a_L+(-\frac{n_L^2d_{L-1}}{1+n_Ld_{L-1}})=b_L+r_L$$
So, $a_L=b_L$ and
$$ r_L=(-\frac{n_L^2d_{L-1}}{1+n_Ld_{L-1}})$$
It is obtained that
$$\frac{1}{r_L}=-\frac{1+n_Ld_{L-1}}{n_L^2d_{L-1}}=b_{L+1}+r_{L+1}$$
Therefore
$$b_{L+1}=\left[-\frac{1+n_Ld_{L-1}}{n_L^2d_{L-1}}\right]$$
Solve the following inequality
\begin{eqnarray*}
&&-\frac{1+n_Ld_{L-1}}{n_L^2d_{L-1}}>K\Leftrightarrow
1+n_Ld_{L-1}>n_L^2(-d_{L-1})K\\
&&\Leftrightarrow (-d_{L-1})<\frac{1}{Kn_L^2+n_L}
\end{eqnarray*}
It is obvious that $(2K+2)n_L(n_L-1)>Kn_L^2+n_L$ for $n_L\ge 2$.
So we take
$$-d_{L-1}<\frac{1}{(2K+2)n_L(n_L-1)}$$
The proof of theorem \ref{theo:basis_alg} is finished.

For practical purpose, we hope the restriction on $n_1>1$ and
$n_1>n_2$ can be lifted. So we have following theorem:
\begin{theorem}\label{theo:practice_alg}
Let $n_0/n_1$ be a reduced rational number and $r$ its
approximation. Assume that $n_0$,$n_1$ are  positive integers and
$N\ge\max\{n_1,2\}$. $K$ is a positive integer. The continued
fraction representations of $n_0/n_1$ and $r$ are
$[a_0,a_1,\cdots,a_L]$ and $[b_0,b_1,\cdots,b_M]$ respectively. If
$|d|=|r-n_0/n_1|<1/((2K+2)N(N-1))$, then one of the following two
statements must hold
\begin{itemize}
\item $a_i=b_i$ for $i=0,\cdots,L$, and $b_{L+1}\ge K$; \item
$a_i=b_i$ for $i=0,\cdots,L-1$, and $b_L=a_L-1$, $b_{L+1}=1$,
$b_{L+2}\ge K$.
\end{itemize}
\end{theorem}
{\bf Proof}: We prove the theorem in three cases:\\ Case 1
($n_1>1,\;n_0<n_1$): This is theorem \ref{theo:basis_alg}.\\
Case 2($n_1=1$): We have that $a_0=n_0/n_1$. If $d=r-n_0/n_1>0$,
then $b_0=a_0$ and $r_0<1/((2K+2)2(2-1))$. So, it holds that
$1/r_0>(2(2K+2))$. Therefore, it is obtained that $b_1>K$. If
$d=r-n_0/n_1<0$, then
\begin{eqnarray*}
&&r=a_0-|d|=a_0-1+1-|d|=b_0+1-|d|\Rightarrow r_0=1-|d|\\
&&\Rightarrow 1/r_0=1+\frac{|d|}{1-|d|}=b_1+\frac{|d|}{1-|d|}\\
&&\Rightarrow
1/r_1=\frac{1-|d|}{|d|}=\frac{1}{|d|}-1>(4K+4)-1\\
&&\Rightarrow b_2>1/r_1-1>4K+2>K
\end{eqnarray*}
So, we have that $b_0=a_0-1$, $b_1=1$ and $b_2>K$.
\\
Case 3($n_0>n_1$): From $n_0/n_1=a_0+n_2/n_1$, it holds that
$n_0/n_1-a_0=n_2/n_1$. On the other hand, we have that
$|n_0/n_1-r|<1/((2K+2)N(N-1))\le 1/n_1$. So, we can deduce that
$a_0<r<a_0+1$. Accordingly, it holds that $b_0=a_0$. Hence, we
have
\begin{eqnarray*}
|d|&=&|r-n_0/n_1|=|b_0+r_0-a_0-n_2/n_1|=|r_0-n_2/n_1|\\
&=&d_0<1/((2K+2)N(N-1))<1/((2K+2)n_1(n_1-1))
\end{eqnarray*}
Since $n_1>1$ and $n_2<n_1$, from theorem \ref{theo:basis_alg},
the theorem holds. Therefore, the proof is finished.

From theorem \ref{theo:practice_alg}, we can get exact
non-negative number $n_2/n_ 1$ from its approximation $r$ as
follows:
\begin{list}{Step \arabic{num}:}{\usecounter{num}\setlength{\rightmargin}{\leftmargin}}
\item estimating an upper bound of the denominator of $n_2/n_1$,
Denoted by $N$; \item computing $$d=\frac{1}{(2N+2)N(N-1))}$$
\item obtaining $r$ by approximate method such that
$|r-n_2/n_1|<d$; \item taking $\varepsilon=1/N$ in algorithm
\ref{algor:1} and calling algorithm \ref{algor:1} to get $b$. So
$b=n_2/n_1$.
\end{list}

\section{Experimental results}
The following examples  run in the platform of Maple 10 and PIV
3.0G, 512M RAM. They take little time for obtaining exact rational
numbers from their approximations, so we do not show time.

\textbf{Example 1.} Let $a$ be unknown rational number. We only
know a bound of its denominator $N=170$. According to theorem
\ref{theo:practice_alg}, Computing rational number $a$ as follows:
Take $K=170+1$, $d=(2*K+2)*N*(N-1)=9883120$, and compute
$1/d=1/9883120$. Assume that we use some numerical method to get
an approximation $b=0.8106507864$ such that  $|a-b|<1/d$. Taking
$\varepsilon=1/K$, we  recover number $a$ by algorithm
\ref{algor:1}. We get $[0,1,4,3,1,1,4,]$ by the first 7 steps.
When doing at step 8, we get $314$, which is larger that $K$. We
stop and return $a=[0,1,4,3,1,1,4]$. It is $137/169$.

\textbf{Example 2.}Let $a$ be unknown rational number. We only
know a bound of its denominator $N=1790$. According to theorem
\ref{theo:practice_alg}, Computing rational number $a$ as follows:
Take $K=1790+1$, yields $d=(2*K+2)*N*(N-1)=11477079040$. Assume
that we use some numerical method to get an approximation
$b=0.1788708777$ such that  $|a-b|<1/d$. Taking $\varepsilon=1/K$,
we  recover number $a$ by algorithm \ref{algor:1}. We get
$[0,5,1,1,2,3,1,6,2,]$ by the first 8 steps. When doing at step 9,
we get $2722$, which is larger that $K$. We stop and return
$a=[0,5,1,1,2,3,1,6,2]$. It is $320/1789$.

\textbf{Example 3.}Let $a$ be unknown rational number. We only
know a bound of its denominator $N=18$. According to theorem
\ref{theo:practice_alg}, Computing rational number $a$ as follows:
Take $K=N+1$, yields $d=(2*K+2)*N*(N-1)=12240$. Assume that we use
some numerical method to get an approximation $1.882434634$ such
that $|a-b|<1/d$. Taking $\varepsilon=1/K$, we  recover number $a$
by algorithm \ref{algor:1}. We get $[1,1,7,1,1]$ by the first 5
steps. When doing at step 6, we get $41$, which is larger that
$K$. We stop and return $a=[1,1,7,1,1]$. It is $32/17$.

\textbf{Example 4.}  This example is an application in obtaining
exact factors from their approximations. Let
$p=-16-56*y-48*z+64*x^2-32*x*y+48*x*z-45*y^2-96*y*z-27*z^2$ be a
polynomial. We want to use approximate method to get its exact
factors over rational number field. First, we transform $p$ to a
monic polynomial as follows:
$$ p=x^2-\frac{1}{2}xy+\frac{3}{4}xz-\frac{45}{64}y^2-\frac{3}{2}yz-\frac{27}{64}z^2-\frac{7}{8}y-\frac{3}{4}z-\frac{1}{4}$$
the least common multiple of denominators of coefficients of
polynomial $p(x,y,z)$ is 64, which is an upper bound\cite{zhang_1}
of denominators of coefficients of the monic factors of polynomial
$p$. taking $K=66$ yields $d=(2*K+2)*65*64=557440$ and
$\varepsilon=1/K$. We use numerical methods to get its approximate
factors as follows\cite{zhang_1}:
$$\bar g_1=1.000000000000x+.6250000000067y+1.124999999530z+.5000000000000$$
$$\bar g_2=1.000000000000x-1.125000000015y-.3749999995480z-.5000000000000$$
the error of coefficients of $\bar g_1$ and $\bar g_2$ is less
than $1/d$ by the numerical methods. According to theorem
\ref{theo:practice_alg}, taking $\varepsilon$ in algorithm
\ref{algor:1}, we obtain two exact factors:
$$g_1=x+\frac{5}{8}y+\frac{9}{8}z+\frac{1}{2}$$
$$g_2=x-\frac{9}{8}y-\frac{3}{8}z-\frac{1}{2}$$

\section{Conclusion}
This paper builds a bridge spanning the gap between approximate
computation and exact results. The exact results can be obtained
by our algorithm as long as we get a bound $N$ of absolute values
of their denominators and their approximations with a error less
than $1/((2N+2)N(N-1))$. Basing on our algorithm, we have succeed
in obtaining exact factors of polynomials from their approximate
factors. Our method can be applied in many aspect, such as proving
inequality statements and equality statements, and computing
resultants, etc. Thus we can take fully advantage of approximate
methods to solve larger scale symbolic computation problems.

\end{document}